\documentclass[leqno,a4paper,12pt]{article}
\usepackage{amsfonts,amssymb}
\usepackage[mathcal,mathbf]{euler}

\def\Z{\mathbb{Z}}
\def\Q{\mathbb{Q}}
\def\R{\mathbb{R}}
\def\C{\mathbb{C}}
\def\F{\mathbb{F}}
\def\Gm{\mathbb{G}_{\mathrm{m}}}

\def\A{\mathbb{A}}
\def\O{{\cal O}}

\def\til#1{\widetilde{#1}}% tilde
\def\ovl#1{\overline{#1}}% overline
\def\pf{{\indent\textit{Proof.}\ }}
\def\qed{\hfill$\square$}
\def\GL{\mathrm{GL}}
\def\GU{\mathrm{GU}}
\def\SL{\mathrm{SL}}
\def\PGL{\mathrm{PGL}}
\def\M{\mathrm{M}}
\def\Spec{\mathop{\mathrm{Spec}}\nolimits}% spectrum
\def\det{\mathop{\mathrm{det}}\nolimits}% determinant
\def\opp#1{{{#1}}^{\mathrm{op}}}% opposite algebra
\def\model{\mathcal{I}}
\def\submodel{\mathcal{J}}
\def\group{\Gamma}
\def\lattice{L}
\def\I{\mathbf{I}}
\def\unitary{U}
\def\Gal{\mathrm{Gal}}
\def\G{\mathbf{G}}
\def\T{\mathbf{T}}
\begin{document}
\newcounter{state}[section]
\setcounter{state}{0}
\setcounter{section}{0}
\renewcommand{\thestate}{\arabic{section}.\arabic{state}}
\renewcommand{\thesection}{\arabic{section}}
\newcommand{\statenumber}{\refstepcounter{state}
\indent{\bf\thestate}}
\newcommand{\sectionnumber}{\refstepcounter{section}
\indent{\bf \thesection.}}
\begin{center}
\begin{large}
{\sf Arithmetic structure of CMSZ fake projective planes}
\footnote{{\em $1991$ Mathematics Subject Classification\/}.
Primary 14G35; Secondary 14J29, 11G15, 11Q25.}
\end{large}

\vspace{4ex}
{\sc Fumiharu Kato}\ and\ {\sc Hiroyuki Ochiai}

\vspace{4ex}
\begin{minipage}{10.5cm}
\setlength{\baselineskip}{.85\baselineskip}
\begin{small}
We prove that the fake projective planes constructed from the diadic
discrete group discovered by Cartwright, Mantero, Steger, and Zappa are
connected Shimura varieties associated to a certain unitary group.
The necessary Shimura data, as well as the unitary group, are explicitly 
described. We also give a field of definition of these fake projective
planes.
\end{small}
\end{minipage}
\end{center}

\sectionnumber\label{section-intro}\ {\bf Introduction}\ 

\vspace{1ex}
In \cite{CMSZ2}, Cartwright, Mantero, Steger, and Zappa discovered a
unitary group in three variables with respect to the quadratic extension 
$\Q(\sqrt{-15})/\Q$ whose integral model over the integer ring with the
prime $2$ inverted gives rise to a diadic discrete group acting
transitively on vertices of Bruhat-Tits building over $\Q_2$.
Inside the integral model are three subgroups to which the
restricted action is simply transitive. 
Moreover, a slight inspection shows that two of them act freely also on
simplicies  
of the other dimensions. Such a situation evokes that Mumford
\cite{Mu79} had also obtained a discrete group with the same properties
(but in a different unitary group)
to construct an algebraic surface with $P_{\mathrm{g}}=q=0$,
$c^2_1=3c_2=9$, and with the ample canonical class, so-called, {\it fake
projective planes}, which is among the most interesting classes of
algebraic surfaces. 
Like that Mumford's group produces such a surface through diadic
uniformization, these two groups give rise to two fake projective
planes. 
It has been shown in \cite{IK98} that these three fake projective planes 
are not isomorphic to among others.

In \cite{Ka99}, in the mean time, it was proved that Mumford's fake
projective plane is a unitary Shimura variety. A clue to this result is
that Mumford's discrete group is, as it is clear from the construction,
a congruence subgroup in the integral model of the unitary group, as well
as that it is arithmetic.  
However, it is not clear whether the groups in \cite{CMSZ2} are
characterized by congruence conditions or not, because the definition of 
the groups therein only gives us generators in matrices, although
they are certainly arithmetic.
This obscurity should be resolved, for it is the key point to settle 
whether the other two fake projective planes also have such a nice
arithmetic structure as Mumford's one has. 

This point is exactly what we will discuss in this paper. In the next
section we will construct congruence subgroups in the unitary group
characterized by conditions in (modulo $3$)-reduction, and will prove
that they coincide up to scalars with the groups defined in \cite{CMSZ2}. 
Properly speaking, we give another way of construction of these groups,
and hence the main stream of our argument is independent from that of
\cite{CMSZ2}, although we have been much motivated by it.
From this we proceed to construct the Shimura varieties in the following 
section, which mimics the argument in \cite{Ka99}.
Our main theorem (Theorem \ref{thm-main}) states that, for each of the
two groups, there exists a Shimura variety with the reflex field
$\Q(\sqrt{-15})$ having two geometrically connected components
isomorphic to the corresponding fake projective plane; moreover, we will 
also see that these connected components have $\Q(\sqrt{-3},\sqrt{5})$
as a field of definition. 

The last section is an appendix to the next section, in which we
will concentrate on proving a technical result.

\vspace{1ex}
{\it Notation and conventions.}\ 
The notation $A\otimes B$ (absence of the base ring) only occurs when 
$A$ and $B$ are $\Q$-algebra, and the tensor is taken over $\Q$.
For a field extension $F/E$, we denote by $\mathrm{Res}_{F/E}$ the Weil
restriction. 
By an involution of a ring $A$ we always mean a homomorphism
$\ast\colon A\rightarrow\opp{A}$ such that $\ast\circ\ast=\mathrm{id}_A$.

\vspace{2ex}
\sectionnumber\label{section-CMSZ}\ 
{\bf The CMSZ-group}  

\vspace{1ex}
\statenumber\label{para-discrete}{\bf .}\ 
{\it Discrete group from unitary group.}\ 
Let $K\subset\C$ be an imaginary quadratic extension of $\Q$, and $p$ a
rational prime subject to the prime decomposition of form
$(p)=\mathfrak{p}\ovl{\mathfrak{p}}$ in $K$ .
We consider the unitary group $\I=\I_Q$ associated to a positive
definite Hermitian matrix $Q$ with entries in $K$; the
group $\I$ is the $\Q$-algebraic group such that 
$$
\I(\Q)=\{g\in\GL_nK\,|\,g^{\ast}Qg=c(g)Q,\ c(g)\in\Q^{\times}\}.
$$
Let $\model$ be the integral model of $\I(\Q)$ consisting of
matrices $Q$-unitary similitudes lying in
$\GL_n\O_K{\textstyle{[\frac{1}{p}]}}$; 
viz., 
$$
\model=\I(\Q)\,{\textstyle{\bigcap}}\,
\GL_n\O_K{\textstyle{[\frac{1}{p}]}}.
$$
Since $p$ decomposes on $K$, the group $\I(\Q_p)$ is, once we choose
$K\hookrightarrow\Q_p$ continuous with respect to the
$\mathfrak{p}$-adic topology on $K$, identified with $\GL_n\Q_p$, and
thus we have a homomorphism $\I(\Q)\rightarrow\PGL_n\Q_p$.
Let $\group$ be the image of $\model$ in $\PGL_n\Q_p$.
Since $\I_{\R}\cong\GU(3)$, as is well-known, $\group$ is discrete and
co-compact in $\PGL_n\Q_p$.

\vspace{1ex}
\statenumber\label{para-CMSZ}{\bf .}\ 
{\it CMSZ-situation.}\ 
The so-called CMSZ-group is the group $\model$ (or $\group$) as above in $n=3$,
$p=2$, and the following setting:
First we set $K=\Q(\sqrt{-15})$.
The integer ring $\O_K$ is $\Z[\lambda]$ with
$\lambda=(1-\sqrt{-15})/2$; the prime $2$ is decomposed as
$(2)=\mathfrak{p}\ovl{\mathfrak{p}}$ where $\mathfrak{p}=2\Z+\lambda\Z$.
The Hermitian matrix $Q$ is the one given by
$$
Q=\left[
\begin{array}{ccc}
\scriptstyle{10}&\scriptstyle{-2(\lambda+2)}&\scriptstyle{\lambda+2}\\
\,\scriptstyle{-2(\ovl{\lambda}+2)}&\scriptstyle{10}&
\,\scriptstyle{-2(\lambda+2)}\\
\scriptstyle{\ovl{\lambda}+2}&\scriptstyle{-2(\ovl{\lambda}+2)}&
\scriptstyle{10}
\end{array}
\right].
$$
The embedding $K\hookrightarrow\Q_2$ is chosen to be the
$\mathfrak{p}$-adic completion; hence, $\lambda/2$ gives a uniformizing
prime, and $\ovl{\lambda}=1-\lambda$ is a diadic unit.
The significance of the group $\group$ thus obtained is the following
fact: 

\vspace{1ex}
\statenumber\label{thm-CMSZ}{\bf .\ Theorem (\cite[\S 5]{CMSZ2}).}\
{\sl The group $\group$ acts transitively on the vertices of the
Bruhat-Tits building $\Delta$ attached to $\PGL_3\Q_2$.}
\qed

\vspace{1ex}
We insert the proof herein for the reader's convenience.
Before doing it, we introduce two elements belonging to
$\model$:
$$
\rho=\left[
\begin{array}{ccc}
\scriptstyle{0}&\scriptstyle{0}&
\frac{\lambda}{2}\\
\scriptstyle{0}&\llap{$\scriptstyle{-}$}\scriptstyle{1}&
{\scriptstyle{1+}}\frac{\lambda}{2}\\
\scriptstyle{1}&\llap{$\scriptstyle{-}$}\scriptstyle{1}&\scriptstyle{1}
\end{array}
\right],\qquad
\tau=\left[
\begin{array}{ccc}
\scriptstyle{0}&\llap{$\scriptstyle{-}$}\scriptstyle{1}&
\frac{\lambda}{2}\\
\scriptstyle{1}&\llap{$\scriptstyle{-}$}\scriptstyle{1}&
{\scriptstyle{1+}}\frac{\lambda}{2}\\
\scriptstyle{0}&\scriptstyle{0}&\scriptstyle{1}
\end{array}
\right].
$$
Note that $\rho^3=\frac{\lambda}{2}I_3$ and $\tau^3=I_3$.

\vspace{1ex}
\pf
Let $\Lambda_0$ be the vertex which is the similarity class of the
standard lattice $(\Z_2)^3$.
Since $\Delta$ is a connected complex, it suffices to show that, for
each vertex $\Lambda$ adjacent to $\Lambda_0$, there exists an element
$g\in\model$ such that $\Lambda=g\Lambda_0$.
There exist elements $g_i\in\model$ for $i\in\Z/7\Z$ with $g_3=\rho$
satisfying the following relations:
$$
\tau^{-1}g_{2i}\tau\ =\ g_i,
\leqno{\textrm{\indent(\ref{thm-CMSZ}.1)}}
$$
$$
\begin{array}{lclclcr}
g_3g_3g_3&=&g_6g_6g_6&=&g_5g_5g_5&=&
\textstyle{\frac{\lambda}{2}}I_3\rlap{,}\\
g_1g_1g_0&=&g_2g_2g_0&=&g_4g_4g_0&=&I_3\rlap{,}\\
g_1g_3g_6&=&g_2g_6g_5&=&g_4g_5g_3&=&I_3\rlap{.}
\end{array}
\leqno{\textrm{\indent(\ref{thm-CMSZ}.2)}}
$$
Then one can easily show that the set $\{g_i^{\pm
1}\Lambda_0\,|\,i\in\Z/7\Z\}$ coincides with the set of vertices
adjacent to $\Lambda_0$. 
\qed

\vspace{1ex}
\statenumber\label{rem-chamber}{\bf .\ Remark (\cite{CMSZ1}\cite{CMSZ2}).}\ 
Let $\mathcal{F}$ be the subset of $(\Z/7\Z)^3$ given by
$$
\mathcal{F}=
\left\{(i,j,k)\in (\Z/7\Z)^3
\,\bigg|\,
\begin{minipage}{4.5cm}
\setlength{\baselineskip}{.85\baselineskip}
\begin{small}
Either $g_ig_jg_k$, $g_jg_kg_i$, or $g_kg_ig_j$ appears in 
(\ref{thm-CMSZ}.2).
\end{small}
\end{minipage}
\right\}.
$$
Then $\mathcal{F}$ consists of $21=3+6\times 3$ elements. Note that 
$g_ig_jg_k$ is either $I_3$ or $\frac{\lambda}{2}I_3$ for each
$(i,j,k)\in\mathcal{F}$.
This set has the following meaning: For each $(i,j,k)\in\mathcal{F}$,
the set $\{g_i\Lambda_0,\Lambda_0,g_k^{-1}\Lambda_0\}$ forms a chamber
in $\Delta$, which we denote by $C(i,j,k)$; this gives rise to a bijection 
between $\mathcal{F}$ and the set of all (non-oriented) chambers
containing $\Lambda_0$.
The following facts are easy to see:

\vspace{1ex}
(\ref{rem-chamber}.1)\ $g_i^{-1}C(i,j,k)=C(j,k,i)$; in particular, $g_3$ 
(resp.\ $g_6$, $g_5$) stabilizes the chamber $C(3,3,3)$ (resp.\
$C(6,6,6)$, $C(5,5,5)$).

\vspace{1ex}
(\ref{rem-chamber}.2)\ $\tau C(i,j,k)=C(2i,2j,2k)$ (this follows from
(\ref{thm-CMSZ}.1) and the fact that $\tau$ fixes $\Lambda_0$).

\vspace{1ex}
Using these facts, one can classify all the finite subgroups in
$\model$ (Theorem \ref{thm-finite} below).

\vspace{1ex}
\statenumber\label{prop-vertex}{\bf .\ Proposition.}\
{\sl Let $\model$ act on $\Delta$ through $\model\rightarrow\PGL_3\Q_2$.
Then the stabilizer $\mathrm{Stab}\Lambda_0$ in $\model$ of the vertex
$\Lambda_0$ is the subgroup generated by $\tau$ and scalar matrices. In
particular, the group $\model$ is generated by $\rho$, $\tau$, and
scalars.} 

\vspace{1ex}
\pf
As the proof of Theorem \ref{thm-CMSZ} indicates, the subgroup of $\model$
generated by $\rho$ and $\tau$ acts transitively on the vertices of
$\Delta$. Hence the second assertion follows from the first. 
Let $g\in\model$ be an element which fixes $\Lambda_0=[(\Z_2)^3]$.
Multiplying by a scalar if necessary, one may assume
$g(\Z_2)^3=(\Z_2)^3$.
There, what we shall prove is:

\vspace{1ex}
\statenumber\label{lem-vertex}{\bf .\ Lemma.}\
{\sl A matrix $g\in\model$ with $g(\Z_2)^3=(\Z_2)^3$ must be of form
$\pm\tau^i$ ($i=0,1,2$).}

\vspace{1ex}
The proof is elementary, but technical, being partly carried out by 
computer search; we postpone it to Appendix. 
\qed

\vspace{1ex}
Recall that a labelling of $\Delta$ is a map $l$ from the set of
vertices to $\Z/3\Z$ which assigns $\nu(\det g)$ mod $3$ to
$\Lambda=[g(\Z_2)^3]$ with $g\in\PGL_3\Q_2$, where $\nu$ is the
normalized valuation $\nu\colon\Q^{\times}_2\rightarrow\Z$.
The map $l$ restricted to each chamber is bijective, and hence $l$ gives 
rise to an orientation in $\Delta$.
Since $l(g\Lambda)=l(\Lambda)+(\nu(\det g)\ \textrm{mod}\ 3)$ for any
$g\in\PGL_3\Q_2$ and any $\Lambda$, the action of $\PGL_3\Q_2$ on
$\Delta$ preserves the orientation.

\vspace{1ex}
\statenumber\label{prop-chamber}{\bf .\ Proposition.}\
{\sl Let us consider, for $(i,j,k)\in\mathcal{F}$, the stabilizer
$\mathrm{Stab}C(i,j,k)$ in $\model$ of the chamber $C(i,j,k)$. Then, if
$i=j=k=3$, $6$, or $5$, $\mathrm{Stab}C(i,j,k)$ is the subgroup
generated by $g_i$ and scalar matrices; otherwise, it consists only of
scalars.} 

\vspace{1ex}
\pf
In view of (\ref{rem-chamber}.1) and (\ref{rem-chamber}.2), it suffices
to show the proposition in the cases $(i,j,k)=(3,3,3)$, $(1,1,0)$, and
$(1,3,6)$. 
Let $g\in\model$ stabilizes $C(i,j,k)$, and assume that $g$ is not a
scalar.
By (\ref{rem-chamber}.2) and Proposition \ref{prop-vertex}, $g$ does not 
fix the vertex $\Lambda_0$, but $g^3$ does, since $g$ preserves the
orientation by the labelling.
Then, taking inverse if necessary, one may assume
$g\Lambda_0=g_i\Lambda_0$.
By Proposition \ref{prop-vertex}, we deduce $g=g_i\tau^jc$ for
some $j\in\{0,1,2\}$ and $c\in(\O_K[\frac{1}{2}])^{\times}$.
Hence, from Proposition \ref{prop-vertex}, it follows that $(g_i\tau^j)^9$
must be a scalar. In case $i=1$, this can be shown to be impossible by
direct calculation. If $i=3$, this is only the case when $j=0$ (hence
$g^3$ is a scalar), as one can check directly.
\qed

\vspace{1ex}
\statenumber\label{thm-finite}{\bf .\ Theorem.}\
{\sl Any non-trivial finite subgroup of $\group$ is conjugate to either
$\langle\rho\rangle$ or to $\langle\tau\rangle$.}

\vspace{1ex}
\pf
Let $H$ be a finite subgroup in $\group$.
Then it is well-known that $H$ stabilizes a simplex in $\Delta$; since
$\group$ preserves the orientation, $H$ stabilizes either a vertex or a
chamber. 
Since $\group$ acts on the vertices transitively, the result follows
from Proposition \ref{prop-vertex} and Proposition \ref{prop-chamber}.
\qed

\vspace{1ex}
\statenumber\label{para-reduction}{\bf .}\ 
{\it Reduction of unitary group.}\ 
For an ideal $\mathfrak{c}$ of $\O_K[\frac{1}{2}]$ which is stable under 
the complex conjugation and a non-negative integer $n$, the Artinian
ring $R_n=\O_K[\frac{1}{2}]/\mathfrak{c}^{n+1}$ has the induced involution,
which we also denote by $\alpha\mapsto\ovl{\alpha}$.
There is the obvious (modulo $\mathfrak{c}^{n+1}$)-reduction morphism
$\pi_n\colon\model\rightarrow\GL_3R_n$, whose image is a subgroup of the 
unitary group 
$$
\unitary_n=\{g\in\GL_3R_n\,|\,g^{\ast}Q_ng=c(g)Q_n,\ c(g)\in F^{\times}_n\},
$$
where $Q_n$ is the (modulo $\mathfrak{c}^n$)-reduction of $Q$ (note that
our $Q$ has coefficients in $\O_K[\frac{1}{2}]$), and $F_n$ is the
subring of $R_n$ consisting of elements fixed by the involution.
Our particular interest is in the case $\mathfrak{c}=3\Z+(\lambda+1)\Z$ 
and $n=0$ or $1$.
Note that $\mathfrak{c}^2=(3)$ gives the prime decomposition of $3$ on
$K$.
Our goal in this section is to construct subgroups in $\model$ which acts
freely on $\Delta$ and transitively on the vertices in terms of congruence
conditions in these reductions. 
To this end, we will first determine the subgroups $\pi_i(\model)$
($i=0,1$) and then, ask for nice subgroups of them of which the inverse
images by $\pi_i$ keep the transitivity on vertices and rule 
out torsion elements.

Before carrying out this program, we change our presentation of matrices
into more convenient form: Let us define a matrix
$$
\Phi=\left[
\begin{array}{ccc}
\scriptstyle{1}&\scriptstyle{2\lambda-1}&
\llap{$\scriptstyle{-}$}\scriptstyle{2}\\
\llap{$\scriptstyle{-}$}\frac{4\lambda-3}{2}&
\frac{2\lambda+1}{2}&
\llap{$\scriptstyle{-}$}\scriptstyle{1}\\
\llap{$\scriptstyle{-}$}\frac{2\lambda-1}{2}&\frac{1}{2}&
\llap{$\scriptstyle{-}$}\scriptstyle{2} 
\end{array}
\right],
$$
and the new Hermitian matrix $Q'=\Phi^{\ast}Q\Phi$, which appears to be:
$$
Q'=\left[
\begin{array}{ccc}
\scriptstyle{90}&\scriptstyle{2\ovl{\lambda}-1}&
\llap{$\scriptstyle{-}$}\scriptstyle{15}\\
\scriptstyle{2\lambda-1}&\scriptstyle{90}&
\scriptstyle{15(2\lambda-1)}\\
\llap{$\scriptstyle{-}$}\scriptstyle{15}&
\scriptstyle{15(2\ovl{\lambda}-1)}&
\scriptstyle{70}
\end{array}
\right].
$$
We set $\model'=\Phi^{-1}\model\Phi$.
(Note that, since $\det\Phi=2^27$, this twist is not admissible over
$\O_K[\frac{1}{2}]$; however, it is harmless because we work only over
the prime above $3$.)
Set $\unitary'_n=\Phi^{-1}\unitary_n\Phi$, and let
$\pi'_n\colon\model\rightarrow\unitary'_n$ be the composite of the
reduction map $\pi_n$ followed by the isomorphism 
$\unitary_n\stackrel{\sim}{\rightarrow}\unitary'_n$. 

\vspace{1ex}
\statenumber\label{para-unitarygroups}{\bf .}\ 
{\it The unitary groups $\unitary'_0$ and $\unitary'_1$.}\ 
First we note that $R_0=F_0\cong\F_3$ and that the induced involution acts
on it trivially. Since
$$
Q'_0=\left[
\begin{array}{ccc}
\scriptstyle{0}&\scriptstyle{0}&\scriptstyle{0}\\
\scriptstyle{0}&\scriptstyle{0}&\scriptstyle{0}\\
\scriptstyle{0}&\scriptstyle{0}&\scriptstyle{1}
\end{array}
\right],
$$
the group $\unitary'_0$ modulo center is isomorphic to the Euclidean
motion group on the affine plane over $\F_3$; more precisely,
$$
\unitary'_0=\left\{
\left[
\begin{array}{cc}
\scriptstyle{A}&\scriptstyle{B}\\
\scriptstyle{0}&\llap{$\scriptstyle{\pm}$}\scriptstyle{1}
\end{array}
\right]\,\bigg|\,
A\in\SL_2\F_3,\ B\in\M_{21}\F_3
\right\}.
\leqno{\textrm{\indent(\ref{para-unitarygroups}.1)}}
$$
The order of $\unitary'_0$ is therefore $432=2^43^3$.

The Artinian ring $R_1$ in turn is isomorphic to $\F_3[t]/(t^2)$, where
$t$ is the image of $1-2\lambda$ ($=\sqrt{-15}$), which is acted on by
the involution $t\mapsto -t$; hence, $F_1\cong\F_3$.
The group $\unitary'_1$ is the unitary group with respect to the
Hermitian form
$$
Q'_1=\left[
\begin{array}{ccc}
&t&\\
-t&\\
&&1
\end{array}
\right],
$$
whence
$$
\unitary'_1=\left\{
\left[
\begin{array}{cc}
\scriptstyle{A}&\scriptstyle{B}\\
\scriptstyle{tC}&\scriptstyle{d}
\end{array}
\right]\,\bigg|\,
\begin{minipage}{8.5cm}
\setlength{\baselineskip}{.85\baselineskip}
\begin{small}
$A\in\GL_2R_1$, $B\in\M_{21}R_1$ and $d\in R_1^{\times}$ 
with $(A$ mod $t)\in\SL_2\F_3$ and ${}^tC=(-d^{-1}JA^{-1}B$ mod 
$t)$
\end{small}
\end{minipage}
\right\},
\leqno{\textrm{\indent(\ref{para-unitarygroups}.2)}}
$$
where $J$ denotes the standard symplectic matrix.
The order of this group is calculated to be $944784=2^43^{10}$.

\vspace{1ex}
\statenumber\label{para-elements}{\bf .}\ 
{\it Special elements.}\ 
Here we introduce some matrices in $\unitary'_1$ which will play
important roles in analyzing the structure of $\unitary'_1$:
$$
z=\left[
\begin{array}{ccc}
1+t\\&1+t\\&&1+t
\end{array}
\right],\ 
u=\left[
\begin{array}{ccc}
1&1\\&1\\&&1-t
\end{array}
\right],\ 
w=\left[
\begin{array}{ccc}
&\llap{$-$}1\\1\\&&1
\end{array}
\right],\ 
$$
$$
b_1=\left[
\begin{array}{ccc}
1&&t\\&1\\&&1
\end{array}
\right],\ 
b_2=\left[
\begin{array}{ccc}
1\\&1&t\\&&1
\end{array}
\right],\ 
c_1=\left[
\begin{array}{ccc}
1&t&1\\&1\\&\llap{$-$}t&1
\end{array}
\right],\ 
c_2=\left[
\begin{array}{ccc}
1\\\llap{$-$}t&1&1\\t&&1
\end{array}
\right],\ 
$$
$$
d_1=\left[
\begin{array}{ccc}
\llap{$1$}+\rlap{$t$}\\&1\\&&1
\end{array}
\right],\ 
d_2=\left[
\begin{array}{ccc}
1\\&\llap{$1$}+\rlap{$t$}\\&&1
\end{array}
\right],\ 
d_3=\left[
\begin{array}{ccc}
1&t\\&1\\&&1
\end{array}
\right],\ 
d_4=\left[
\begin{array}{ccc}
1\\t&1\\&&1
\end{array}
\right].
$$
We set
$$
T\colon=\langle b_1,b_2\rangle,\ 
H\colon=\langle z,c_1,c_2\rangle,\ 
M\colon=\langle d_1,d_2,d_3,d_4\rangle,\ 
S\colon=\langle u,w\rangle.
$$

\statenumber\label{ntn-operations}{\bf .\ 
Notation.}\ We hereafter make use of the following notation:
\begin{eqnarray*}
[\alpha,\beta]&=&\alpha^{-1}\beta^{-1}\alpha\beta,\\
\alpha^{\delta}&=&\delta^{-1}\alpha\delta.
\end{eqnarray*}

\statenumber\label{para-relations}{\bf .}\ 
{\it Relations among special elements.}\ 
Here we list up the relations among those elements defined above:

\vspace{1ex}
(\ref{para-relations}.1)\ All but $w$ are of order $3$, while the order
of $w$ is $4$.

\vspace{1ex}
(\ref{para-relations}.2)\ $z$ lies in the center of $\unitary'_1$, and 
$[c_1,c_2]=z$. Therefore, the subgroup $H$ is the Heisenberg group with
the center $\langle z\rangle$ and $H/\langle z\rangle\cong\F^2_3$.

\vspace{1ex}
(\ref{para-relations}.3)\ 
$T$ and $M$ are vector groups; $T\cong\F^2_3$ and $M\cong\F^4_3$.
We moreover have $[T,H]=[T,M]=1$.

\vspace{1ex}
(\ref{para-relations}.4)\ In $S$ we have the relations
$[u,w^2]=1,\ wuw=u^{-1}wu^{-1},\ wu^{-1}w=uw^{-1}u$.
Hence $S$ is isomorphic to $\SL_2\F_3$.

\vspace{1ex}
(\ref{para-relations}.5)\ $T$ is normalized by $S$; actually, we have 
$b_1^u=b_1,\ b_2^u=b_1^{-1}b_2,\ b_1^w=b_2^{-1},\ b_2^w=b_1$.

\vspace{1ex}
(\ref{para-relations}.6)\ The action of $S$ on $H$ is given by
$c_1^u=b_1^{-1}c_1,\ c_2^u=b_1^{-1}c_1^{-1}c_2,\ c_1^w=c_2^{-1},\ 
c_2^w=c_1$. 
Hence, in particular, the $\langle T,H\rangle$ is normalized by $S$.

\vspace{1ex}
(\ref{para-relations}.7)\ The action of $H$ on $M$ is by
$$
\begin{array}{llll}
d_1^{c_1}=b_1d_1,&d_2^{c_1}=d_2,&d_3^{c_1}=d_3,&d_4^{c_1}=b_2d_4,\\
d_1^{c_2}=d_1,&d_2^{c_2}=b_2d_2,&d_3^{c_2}=b_1d_3,&d_4^{c_2}=d_4.\\
\end{array}
$$

(\ref{para-relations}.8)\ $M$ is normalized by $S$, since
$$
\begin{array}{llll}
d_1^u=d_1d_3,&d_2^u=d_2d_3^{-1},&d_3^u=d_3,&d_4^u=d_1^{-1}d_2d_3^{-1}d_4,\\
d_1^w=d_2,&d_2^w=d_1,&d_3^w=d_4^{-1},&d_4^w=d_3^{-1}.\\
\end{array}
$$

\statenumber\label{prop-unitary_1}{\bf .\ Proposition.}\
{\sl The group $\unitary'_1$ is generated by the subgroups $T$, $H$,
$M$, $S$, and $-I_3$. Moreover, there exists an isomorphism}
$$
\unitary'_1\stackrel{\sim}{\longrightarrow}
(((T\times M)\rtimes H)\rtimes S)\times\F^{\times}_3
\stackrel{\sim}{\longrightarrow}((\F_3^6\rtimes
H)\rtimes\SL_2\F_3)\times\F^{\times}_3.
$$

\pf
Let $\unitary^{\prime +}_1$ be the subgroup generated by $T$, $H$, $M$,
and $S$. 
By the appearance of the matrices in \ref{para-elements} and relations in
\ref{para-relations}, we easily see  
$\unitary^{\prime +}_1\cong((T\times M)\rtimes H)\rtimes S$.
But the order of the left-hand side already attains the half of that of
$\unitary'_1$, that is, $2^33^{10}$.
Hence $\unitary^{\prime +}_1$ is precisely the kernel of the surjective
homomorphism 
$$
(\textrm{det\ mod}\ t)\colon\unitary'_1\longrightarrow\F^{\times}_3.
\leqno{\textrm{\indent(\ref{prop-unitary_1}.1)}}
$$
Clearly, this morphism has a cross-section with the image in the center
of $\unitary'_1$, thereby the assertion. 
\qed

\statenumber\label{prop-model_1}{\bf .\ Proposition.}\
{\sl The image $\pi'_1(\model)$ by the (modulo $3$)-reduction of
$\model$ is the subgroup of $\unitary'_1$ generated by $T$, $H$, $S$,
$d_1d_2$, and scalars (hence is of order $2^43^7=34992$).
In particular, we have $\pi'_0(\model)=\unitary'_0$.}

\vspace{1ex}
\pf By Proposition \ref{prop-vertex}, we know that $\pi'_1(\model)$ is
generated by $(\Phi^{-1}\rho\Phi\ \textrm{mod}\ 3)$,
$(\Phi^{-1}\tau\Phi\ \textrm{mod}\ 3)$ (which we denote, in this proof,
by $\rho$ and $\tau$, respectively, for simplicity), and scalars.
As one checks easily, 
$$
\begin{array}{lcccl}
\rho&=&\left[
\begin{array}{ccc}
\scriptstyle{0}&\scriptstyle{-1+t}&\scriptstyle{-t}\\
\scriptstyle{1-t}&\scriptstyle{-1-t}&\scriptstyle{1+t}\\
\scriptstyle{0}&\scriptstyle{-t}&\scriptstyle{1+t}
\end{array}
\right]&=&b_1^{-1}c_2wu^{-1}(d_1d_2)^{-1},\\
\tau&=&\left[
\begin{array}{ccc}
\scriptstyle{1+t}&\scriptstyle{1-t}&\scriptstyle{1+t}\\
\scriptstyle{0}&\scriptstyle{1+t}&\scriptstyle{0}\\
\scriptstyle{0}&\scriptstyle{-t}&\scriptstyle{1+t}
\end{array}
\right]&=&z^{-1}c_1u(d_1d_2)^{-1},
\end{array}
\leqno{\textrm{\indent(\ref{prop-model_1}.1)}}
$$
and hence it follows that $\pi'_1(\model)\subseteq\langle
T,H,S,d_1d_2\rangle$. 
For the converse, one calculates
$$
z=\rho^{-3},\quad w=\rho^4(\tau\rho^{-1}\tau^{-1}\rho^{-1})^2,
$$
$$
b_1=[\tau\rho\tau\rho^2\tau\rho\tau^{-1},
\rho\tau\rho\tau^{-1}\rho\tau\rho\tau],\quad
b_2=wb_1w^{-1},
$$
$$
c_1=\rho^{-3}(\tau^{-1}\rho\tau\rho\tau^{-1})^{(\rho\tau\rho)^{-1}}
(\tau\rho)^{-2}b_1b_2^{-1},\quad
c_2=wc_1w^{-1},
$$
$$
u=\tau\rho^{-1}\tau^{-1}\rho^{-1}\tau^{-1}b_1^{-1}b_2c_2^{-1}w^2,
$$
whence $\pi'_1(\model)\supseteq\langle
T,H,S,d_1d_2\rangle$. 
\qed

\vspace{1ex}
Now, to find subgroups of $\model$ as in \ref{para-reduction}, one has
to first look for subgroups which does not contain
$\model$-conjugates of $\tau$ and $\rho$; moreover, the indices of
such subgroups in $\model$ are required to be $3$, since they have to act 
on the vertices of $\Delta$ simply transitively.
This leads us to the following kind of statement:

\vspace{1ex}
\statenumber\label{lem-index3}{\bf .\ Lemma.}\
{\sl Every subgroup of $\pi'_1(\model)$ of index $3$ is conjugate to
either one of the following subgroups:
\begin{eqnarray*}
\submodel_1&=&\langle T,H,P,u,\textrm{scalars}\rangle\\
\submodel_2&=&\langle T,H,P,d_1d_2,\textrm{scalars}\rangle\\
\submodel_3&=&\langle T,H,P,ud_1d_2,\textrm{scalars}\rangle\\
\submodel_4&=&\langle T,H,P,u(d_1d_2)^{-1},\textrm{scalars}\rangle
\rlap{,}
\end{eqnarray*}
where $P$ is the unique $2$-Sylow subgroup of $S$, generated by 
$w$ and $w^u$.}

\vspace{1ex}
\pf
Let $\submodel$ be a subgroup of $\pi'_1(\model)$ of index $3$.
Since $\langle P,\textrm{scalars}\rangle\cong P\times\F^{\times}_3$ is a 
$2$-Sylow subgroup of $\pi'_1(\model)$, replacing it by a suitable
conjugate, we may assume that $\submodel$ contains $\langle
P,\textrm{scalars}\rangle$.
Let us consider the (modulo $t$)-reduction map
$\psi\colon\pi'_1(\model)\rightarrow\pi'_0(\model)$, of which the kernel 
is $K=\langle T,z,d_1d_2\rangle$.
Since $\submodel\bigcap K$ is a normal subgroup in $\submodel$, it is in 
particular a $(P\times\F^{\times}_3)$-module by conjugation.
Since the $(P\times\F^{\times}_3)$-module $K$ is decomposed into
irreducible factors as $T\times\langle z\rangle\times\langle
d_1d_2\rangle\cong\F^2_3\times\F_3\times\F_3$, we deduce that
$\submodel\bigcap K$ contains $T$;
$$
T\subset \submodel{\textstyle \bigcap} K\subset\submodel.
\leqno{\textrm{\indent(\ref{lem-index3}.1)}}
$$

Next we claim that 
$$
\psi(H)\subset\psi(\submodel).
\leqno{\textrm{\indent(\ref{lem-index3}.2)}}
$$
Consider the projection $\pi'_0(\model)\rightarrow
\SL_2\F_3\times\F^{\times}_3$ which extracts the 
top-left $2\times 2$ component and the scalars.
Since its kernel is $\psi(H)$ consisting of $9$ elements, and since
$\psi(\submodel)$ is of index at most $3$ in $\pi'_0(\model)$, we deduce
$\psi(\submodel)\bigcap\psi(H)\neq\{1\}$.
This intersection is, moreover, stable under $P$-conjugation.
But, as one sees easily, the group $\psi(H)$ is an irreducible
$P$-module, thereby the claim.

By (\ref{lem-index3}.2), we know that, for each $i=1,2$, there exists
$k_i\in K$ such that $k_ic_i\in\submodel$.
Since $K$ is an abelian group, we deduce that
$z=[k_1,k_2][c_1,c_2]=[k_1c_1,k_2c_2]\in\submodel$. 
Combining with (\ref{lem-index3}.1), we get $\langle z,T\rangle\subset
\submodel\bigcap K\subset\submodel$.
By this and that $K=\langle T,z,d_1d_2\rangle$, we may take $k_i$ 
from $\langle d_1d_2\rangle$.
Since $[d_1d_2,w^2]=1$, we get $c_1=[c_1,w^2]=[k_1c_1,w^2]\in\submodel$.
Then we moreover have $c_2=wc_1w^{-1}\in\submodel$.
This implies that $\submodel$ contains also $H$, and hence:
$$
\langle T,H,P,\textrm{scalars}\rangle\subset \submodel{\textstyle
\bigcap} K\subset\submodel. 
\leqno{\textrm{\indent(\ref{lem-index3}.3)}}
$$
The left-hand side of (\ref{lem-index3}.3) is already a subgroup in 
$\pi'_1(\model)$ of index $9$.
Then it is now obvious that any subgroup of $\pi'_1(\model)$ of index
$3$ which contains it is one of $\submodel_i$ (i=1,2,3,4)
as stated.
\qed

\vspace{1ex}
\statenumber\label{cons-groups}{\bf .\ Construction.}\
Now, let us define the subgroup $\model_i$ ($i=1,2,3,4$) of $\model$ by 
$$
\model_i=(\pi'_1)^{-1}(\submodel_i),
$$
and let $\group_i$ be its image in $\PGL_3\Q_2$.
These are subgroups of index $3$.
By (\ref{prop-model_1}.1) we know that the subgroup $\group_3$ contains 
$\rho$ and that $\group_4$ contains $\tau$; hence these are not
torsion-free.
The groups $\group_1$ and $\group_2$, on the other hand, has no element
of finite order except $1$, since, as one sees from the relations in
\ref{para-relations}, they does not have any conjugates of $\tau$ and
$\rho$.
Hence the discrete groups $\group_1$ and $\group_2$ act freely on
$\Delta$ and simply transitively on the vertices.

Let $\Omega$ be the Drinfeld symmetric space over $\Q_2$ of dimension $2$.
By \cite[1]{Mu79}, we deduce that both $\group_1\backslash\Omega$ and
$\group_2\backslash\Omega$ are, respectively, algebraized to fake
projective planes ${}^IX_{\mathrm{CMSZ}}$ and ${}^{II}X_{\mathrm{CMSZ}}$
over $\Q_2$, not isomorphic to each other, since our definition shows
that $\group_1$ and $\group_2$ are not conjugate in $\PGL_3\Q_2$ (cf.\
\cite{IK98}). 

\vspace{1ex}
\statenumber\label{rem-CMSZ}{\bf .\ Remark.}\ 
One can check by direct calculation that our groups $\group_i$ for $i=1,2,3$
are exactly the $\Phi$-conjugate of the groups $\group_{B.i}$ 
discovered by Cartwright, Mantero,
Steger, and Zappa in their paper \cite[pp.\ 182]{CMSZ2}.
Indeed, the elements $\rho$ and $\tau$ in \ref{thm-CMSZ} are 
exactly $g_3s$ and $s$ therein, respectively.
The set $\mathcal{F}$ defined in Remark \ref{rem-chamber}, which was
used to analyze the action of $\model$ on the chambers of $\Delta$, is
nothing but the triangle presentation $\mathcal{F}_{B.3}$.

\vspace{2ex}
\sectionnumber\label{section-shimura}\ 
{\bf The Shimura varieties}  

\vspace{1ex}
\statenumber\label{para-situation}{\bf .}\ 
{\it Situation.}\ 
In this section we regard the field $K=\Q(\lambda)$ as a subfield of the
45th-cyclotomic field $\Q(\zeta_{45})$ which is embedded in $\C$ so that
$\zeta_{45}^{21}+\zeta_{45}^{33}+\zeta_{45}^{39}+\zeta_{45}^{42}$ is mapped
to $\lambda$ (cf.\ \ref{para-CMSZ}).
Consider the intermediate field $\Q(\zeta_9)$ with
$\zeta_9=\zeta_{45}^5$, and let $L$ be the composite field of $K$ and
the maximal real subfield of 
$\Q(\zeta_9)$; $L$ is an extension of $K$ of degree $3$ generated $K$ by
$\zeta_9+\zeta_9^{-1}$.
The Galois group $\Gal(L/K)$ is generated by
$\sigma\colon\zeta_9+\zeta_9^{-1}\mapsto\zeta_9^2+\zeta_9^{-2}$.
We need to construct

\vspace{1ex}
$\bullet$ a central division algebra $D$ over $K$,

$\bullet$ a positive involution $\ast$ of $D$ of second kind,

$\bullet$ a non-degenerate anti-Hermitian form $\psi$,

\vspace{1ex}\noindent
from which we will construct the associated Shimura variety.
The construction of these data is basically parallel to that in
\cite[2]{Ka99}, and here, we sometimes limit ourselves to sketchy
presentation.

\vspace{1ex}
\statenumber\label{para-algebra}{\bf .}\ 
{\it Division algebra with involutions.}\ 
First we let 
$$
\mu=\lambda/\ovl{\lambda}\in K.
$$
The division algebra $D$ is defined by
$$
D=\bigoplus^2_{i=0}L\Pi^i;\quad
\Pi^3=\mu,\ \Pi z=z^{\sigma}\Pi\ \textrm{for $z\in L$}.
$$
Then $D\otimes_KL$ is isomorphic to the matrix algebra $M_3(L)$; more
explicitly, such an isomorphism is given by the morphism $D\rightarrow
M_3(L)$ defined by
$$
z\mapsto
\left[
\begin{array}{ccc}
z\\
&z^{\sigma^2}\\
&&z^{\sigma}
\end{array}
\right]\quad
\textrm{and}\quad
\Pi\mapsto
\left[
\begin{array}{ccc}
&&\mu\\
1\\
&1
\end{array}
\right].
\leqno{\textrm{\indent(\ref{para-algebra}.1)}}
$$
By this, $D$ is regarded as a subalgebra over $K$ of $M_3(L)$.
The Hasse invariants of $D$ are zero except for those at $\mathfrak{p}$
and $\ovl{\mathfrak{p}}$, which are $1/3$ and $-1/3$, respectively.
The subring 
$$
\O_D=\O_L\oplus\O_L\ovl{\lambda}\Pi\oplus\O_L\ovl{\lambda}\Pi^2
$$
is an order of $D$, maximal both at $\mathfrak{p}$ and
$\ovl{\mathfrak{p}}$. 
The involution $\ast$ on $D$ is defined by
$$
\Pi^{\ast}=\ovl{\mu}\Pi^2\quad
\textrm{and}\quad
z^{\ast}=\ovl{z}\ \textrm{for}\ z\in L, 
$$
and the anti-Hermitian form $\psi$ by
$\psi(\alpha,\beta)=\mathrm{tr}_{D/\Q}(\alpha^{\ast}b\beta)$, where
$$
b=(\ovl{\lambda}-\lambda)-\ovl{\lambda}\Pi+\ovl{\lambda}\Pi^2.
$$
Since $\mathrm{Nm}^{\circ}_{D|\Q}(b)=-7(\ovl{\lambda}-\lambda)$ is both
$\mathfrak{p}$-adic and $\ovl{\mathfrak{p}}$-adic unit, the pairing
$\psi$ restricted to $\O_D$ is perfect.
The element $b$ also induces another involution $\bigstar$ by
$\alpha^{\bigstar}=b^{-1}\alpha^{\ast}b$.

\vspace{1ex}
\statenumber\label{para-shimura}{\bf .}\ 
{\it Shimura variety.}\ 
Let us consider the $\Q$-algebraic group $\G$ such that
$$
\G(\Q)=\{\gamma\in D^{\times}\,|\,\gamma^{\bigstar}\gamma\in\Q^{\times}\}.
$$
By the identification (\ref{para-algebra}.1) we easily see that
$\G_{\R}$ is isomorphic to $\GU(2,1)$.
The Shimura data is given by the homomorphism 
$$
h\colon\mathrm{Res}_{\C/\R}{\Gm}_{,\C}\longrightarrow\G_{\R}
$$
which sends $\sqrt{-1}$ to
$\mathrm{diag}(\sqrt{-1},\sqrt{-1},-\sqrt{-1})$, and is identity on
scalars.
This gives rise to a Shimura variety
$$
\mathcal{S}h_C=\G(\Q)\backslash X_{\infty}\times\G(\A_{\mathrm{f}})
/C,
$$
where $X_{\infty}$ is the set of all conjugates of $h$, and $C$ is a
(sufficently small) open compact subgroup of $\G(\A_{\mathrm{f}})$.
This is a finite disjoint union of quasi-projective manifolds, which are 
arithmetic quotients of the complex unit-ball.
The reflex field $E$ of $\mathcal{S}h$ is $K$.

\vspace{1ex}
\statenumber\label{pro-inner}{\bf .\ Proposition.}\ 
{\sl Let $\I$ and $\G$ be the $\Q$-algebraic groups defined in
\ref{para-CMSZ} and \ref{para-shimura}, respectively.
Then we have the following:

(a)\ $\I$ is an inner form of $\G$.

(b)\ $\I(\A^2_{\mathrm{f}})\cong\G(\A^2_{\mathrm{f}})$.

\noindent
Moreover, the algebraic group $\I$ is uniquely determined up to
$\Q$-isomorphisms by (a), (b) and the following condition:

(c)\ $\I(\Q_2)\cong\{(x,y)\in\GL_3\Q_2\times
\opp{(\GL_3\Q_2)}\,|\,xy\in\Q_2\}$. }

\vspace{1ex}
\pf
The proof is, literally, almost the same as that of \cite[4.1]{Ka99};
it differs only when we compare, in the proof of (b), two involutions
induced from $Q$ and $b$. We should hereupon compare two Hermitian
matrices $3Q$ and $-7b$ (where $b$ is regarded as a matrix by
(\ref{para-algebra}.1)). The determinants of these are square
numbers, thereby the assertion.
\qed

\vspace{1ex}
\statenumber\label{para-level}{\bf .}\ 
{\it CMSZ-level structures.}\ 
We fix an isomorphism
$\phi\colon\I(\A^2_{\mathrm{f}})\cong\G(\A^2_{\mathrm{f}})$. 
Let $C^{\mathrm{max}}_2$ be the maximal open compact subgroup of
$\G(\Q_2)$ consisting of elements which map $\O_D\otimes_{\Z}\Z_2$ into 
itself.
For a finite set of rational primes $W$ with $2\not\in W$ (resp.\ a
rational prime $p\neq 
2$), we define $C^W_{\mathrm{max}}$ (resp.\ $C^{\mathrm{max}}_p$) to be
the maximal open compact subgroup of $\I(\A^W_{\mathrm{f}})$ (resp.\
$\I(\Q_p)$) consisting of elements which map
$\O_K^3\otimes_{\Z}\widehat{\Z}^W$ (resp.\ $\O_K^3\otimes_{\Z}\Z_p$)
into itself, where $\widehat{\Z}^W$ is the product of all $\Z_{\ell}$'s
for $\ell$ not belonging to $W$.
We look at $C^{\mathrm{max}}_3$, which has the reduction map (also
denoted by $\pi'_1$) to the finite group $\unitary'_1$ (cf.\
\ref{cons-groups}). 
Set ${}^IC_3=(\pi'_1)^{-1}(\submodel_1)$ and
${}^{II}C_3=(\pi'_1)^{-1}(\submodel_2)$, where the groups $\submodel_i$
are those defined in \ref{lem-index3}.
Now, we define
$$
{}^IC^2_{\mathrm{CMSZ}}=C^{2,3}_{\mathrm{max}}{}^IC_3\quad
\textrm{and}\quad
{}^{II}C^2_{\mathrm{CMSZ}}=C^{2,3}_{\mathrm{max}}{}^{II}C_3.
$$
Obviously, the image of
$\I(\Q)\bigcap(\I(\Q_2)\times{}^IC^2_{\mathrm{CMSZ}})$ 
(resp.\ $\I(\Q)\bigcap(\I(\Q_2)\times{}^{II}C^2_{\mathrm{CMSZ}})$) in 
$\I_{\mathrm{ad}}(\Q_2)\cong\PGL_3\Q_2$ is nothing but $\group_1$
(resp.\ $\group_2$).
We finally set
$$
{}^IC=C^{\mathrm{max}}_2\phi({}^IC^2_{\mathrm{CMSZ}})\quad
\textrm{and}\quad
{}^{II}C=C^{\mathrm{max}}_2\phi({}^{II}C^2_{\mathrm{CMSZ}}),
$$
which are open compact subgroups of $\G(\A_{\mathrm{f}})$.

\vspace{1ex}
\statenumber\label{thm-main}{\bf .\ Theorem.}\
{\sl The canonical model $\mathrm{Sh}_{{}^IC}$ (resp.\
$\mathrm{Sh}_{{}^{II}C}$) of the Shimura variety $\mathcal{S}h_{{}^IC}$
(resp.\ $\mathcal{S}h_{{}^{II}C}$) consists of two connected components
defined over $E'=\Q(\sqrt{-3},\sqrt{5})$ which are permuted by the
Galois action of $\Gal(E'/E)$.
Moreover, the base change of each connected component to
$E'_{\mathfrak{p}}\cong\Q_{2^2}$ is isomorphic to the fake projective
plane ${}^IX_{\mathrm{CMSZ}}\times_{\Spec\Q_2}\Spec\Q_{2^2}$
(resp.\ ${}^{II}X_{\mathrm{CMSZ}}\times_{\Spec\Q_2}\Spec\Q_{2^2}$)
defined in \ref{cons-groups}.}

\vspace{1ex}
\pf
The proof is parallel to that of \cite[4.3]{Ka99}; the main part is to
apply \cite[6.5]{RZ96} to our situation.
The only point which should be mentioned here is the computation of the
morphism $\vartheta\colon\G\rightarrow\T$ by
$\gamma\mapsto(\mathrm{Nm}^{\circ}_{\opp{D}|K}(\gamma),c(\gamma))$, where
$$
\T=\{(k,f)\in\mathrm{Res}_{K/\Q}{\Gm}_{,K}\times{\Gm}_{,\Q}\,|\,k\ovl{k}=f^3\}
\cong\mathrm{Res}_{K/\Q}{\Gm}_{,K}.
$$
We shall claim that the images $\vartheta({}^IC)$ and
$\vartheta({}^{II}C)$ are maximal open compact subgroups of
$\mathrm{Res}_{K/\Q}{\Gm}_{,K}(\A_{\mathrm{f}})$; once we prove this,
the assertions on connected components and the field of definition
follow from the fact that $E'=\Q(\sqrt{-3},\sqrt{5})$ is the Hilbert
class field of $E=K$.
To this end, we only have to look at the component at $3$.
What to prove is that the homomorphism $\vartheta\circ\phi$ localized at $3$
$$
\vartheta\circ\phi\colon
{}^IC_3\ \textrm{or}\ {}^{II}C_3\ni k\theta\mapsto
(k\ovl{k})^{-1}k^3\det\theta\in\til{\Z}^{\times}_3,
$$
where $\til{\Z}_3$ is the integer ring of the quadratic ramified
extension of $\Q_3$, is surjective (here $\theta$ has been taken such that
$\theta^{\ast}Q\theta=Q$).

Here we prove it for the group ${}^IC_3$; for the other one is similar, 
even easier.
Since $\Z^{\times}_3$ is contained in ${}^IC_3$, we have 
$\Z^{\times}_3\subset\vartheta\circ\phi({}^IC_3)$.
Also, since $\til{\Z}^{\times}_3$ is contained in ${}^IC_3$, we know
the elememts of form $k^3\ovl{k}$ for $k\in\til{\Z}^{\times}_3$ belongs
to $\vartheta\circ\phi({}^IC_3)$; but since $k\ovl{k}\in\Z^{\times}_3$
is in it, we deduce
$\{k^3\,|\,k\in\til{\Z}^{\times}_3\}
=\pm 1+3\til{\Z}_3\subset\vartheta\circ\phi({}^IC_3)$.
Hence it suffices to prove that $(\vartheta\circ\phi({}^IC_3)$ mod $3)$
is the whole $(\F_3[t]/(t^2))^{\times}\cong(\til{\Z}_3/(3))^{\times}$.
But this can be checked easily; for instance, factoring through
$\unitary'_1$, we get the value $1-t$ by the element $u$ (cf.\
\ref{para-elements}), and hence get $1+t$ and $1$. Multiplying with
$-1$, we get the all.
\qed

\vspace{1ex}
\statenumber\label{cor-main}{\bf .\ Corollary.}\
{\sl Each connected component of the complex surface
$\mathcal{S}h_{{}^IC}$ and $\mathcal{S}h_{{}^{II}C}$ is a fake
projective plane, having a field of definition $\Q(\sqrt{-3},\sqrt{5})$.
Moreover, it is an arithmetic quotient of the complex unit-ball.}
\qed

\vspace{2ex}
\sectionnumber\label{section-appendix}\ 
{\bf Appendix: Proof of Lemma \ref{lem-vertex}}  

\vspace{1ex}
\statenumber\label{step-1}{\bf .}\ 
Let $G_1=\{g\in\model\,|\,g(\Z_2)^3=(\Z_2)^3\}$ and
$G_2=\{\pm\tau^i\,|\,i=0,1,2\}$. What to prove is the equality
$G_1=G_2$. The inclusion $G_2\subset G_1$ is clear.
Let $g$ be an element in $G_1$.
Then $g$ has entries in
$$
\O_K[{\textstyle \frac{1}{2}}]\bigcap\Z_2=
\Z[\lambda,{\textstyle \frac{1}{2}}]\bigcap\Z_2=
\Z[\ovl{\lambda},\ovl{\lambda}^{-1}].
$$
Since $g$ is $Q$-unitary, and since $30Q^{-1}$ takes its entries in
$\Z[\ovl{\lambda},\ovl{\lambda}^{-1}]$, they actually belongs to
$\Z[\ovl{\lambda},\ovl{\lambda}^{-1}]\bigcap
\frac{1}{30}\Z[\lambda,\lambda^{-1}]=
\Z+\frac{\lambda}{2}\Z$, which we denote by $\lattice$.
Set
$$
V=\left\{
v\in\lattice^3\,\bigg|\,
\begin{minipage}{4cm}
\setlength{\baselineskip}{.85\baselineskip}
\begin{small}
$v^{\ast}Qv=10$, $\tau v$ and $\tau^{-1} v$ lie in $\lattice^3$, 
and $v\not\in\frac{\lambda}{2}\Z^3$
\end{small}
\end{minipage}
\right\}.
$$
Then each column vetcor of $g$ belongs to $V$; indeed, $v^{\ast}Qv=10$
and $\tau v$, $\tau^{-1} v\in\lattice^3$ are obvious, and, if it were in 
$(\lambda/2)\Z^3$, then $\det g$ would be a multiple of 
$\lambda/2$, which is not invertible in $\Z_2$.

\vspace{1ex}
\statenumber\label{cla-24vectors}{\bf .\ Claim.}\ 
{\sl The set $V$ consists of $24$ vectors among
$$
G_2\left[
\begin{array}{c}
\frac{\lambda}{2}\\
\scriptstyle{1}\\
\scriptstyle{0}
\end{array}
\right],\quad
G_2\left[
\begin{array}{c}
\scriptstyle{0}\\
\scriptstyle{0}\\
\scriptstyle{1}
\end{array}
\right],\quad
G_2\left[
\begin{array}{c}
\scriptstyle{0}\\
\scriptstyle{1}\\
\scriptstyle{1}
\end{array}
\right],\quad\textrm{and}\quad
G_2\left[
\begin{array}{c}
\llap{$\scriptstyle{-}$}\scriptstyle{1}\\
\llap{$\scriptstyle{-}$}\scriptstyle{1}\\
\scriptstyle{0}
\end{array}
\right].
$$
}

The proof of the claim is done in the following several paragraphs: 

\vspace{1ex}
\statenumber\label{step-2}{\bf .}\ 
For $x=a+(\lambda/2)b$ ($a,b\in\Z$) we have
$x\ovl{x}=a^2+(1/2)ab+b^2\geq(3/4)(a^2+b^2)$.
Let $v={}^t(x_1,x_2,x_3)\in\lattice^3$ with
$x_i=a_i+(\lambda/2)b_i$ ($i=1,2,3$).
Since the characteristic polynomial of $Q$ is $t^3-30t^3+210t-300$, one
sees that the minimum eigenvalue of $Q$ is greater than $1.92$; therefore,
\begin{eqnarray*}
10\ =\ v^{\ast}Qv&\geq&1.92|v|^2\\
&\geq&1.92\cdot(3/4)\cdot({\textstyle \sum}(a_i^2+b_i^2)),
\end{eqnarray*}
whence 
$$
{\textstyle \sum}(a_i^2+b_i^2)\leq 6.
\leqno{\textrm{\indent(\ref{step-2}.1)}}
$$
In particular, we deduce $|a_i|$, $|b_i|\leq 2$.
We are to strengthen this restraint into $|a_i|$, $|b_i|\leq 1$.

\vspace{1ex}
\statenumber\label{step-3}{\bf .}\ 
To this end, let us consider the $\Q$-quadratic form
$F(a_1,a_2,a_3,b_1,b_2,b_3)=v^{\ast}Qv$; more explicitly, 
\begin{eqnarray*}
{\scriptstyle{\frac{1}{10}}}F&=&a_1^2+a_2^2+a_3^2+b_1^2+b_2^2+b_3^2\\
&&-a_1a_2-a_2a_3+{\scriptstyle{\frac{1}{2}}}a_3a_1-b_1b_2-b_2b_3+
{\scriptstyle{\frac{1}{2}}}b_3b_1\\
&&+{\scriptstyle{\frac{1}{2}}}a_1b_1+{\scriptstyle{\frac{1}{2}}}a_1b_2
-{\scriptstyle{\frac{1}{4}}}a_1b_3-a_2b_1\\
&&+{\scriptstyle{\frac{1}{2}}}a_2b_2+{\scriptstyle{\frac{1}{2}}}a_2b_3
+{\scriptstyle{\frac{1}{2}}}a_3b_1-a_3b_2+{\scriptstyle{\frac{1}{2}}}a_3b_3.
\end{eqnarray*}
Then we see, writing
$(a_1,a_2,a_3,b_1,b_2,b_3)=(y_1,y_2,y_3,y_4,y_5,y_6)$, 
$$
{\scriptstyle{\frac{1}{10}}}F\geq
{\textstyle \sum}y_i^2-{\textstyle \sum_{i<j}}|y_iy_j|,
$$
where the right-hand side (which we denote by $F'$) is symmetric in
$y_i$'s.
Assume now that one of $a_i$'s or $b_i$'s has absolute value $2$. 
Due to (\ref{step-2}.1) and the symmetry, the possible values of $F'$
are $F'(2,0,0,0,0,0)=4$, $F'(2,1,0,0,0,0)=3$, and $F'(2,1,1,0,0,0)=1$.
Hence, in order that $(1/10)F=1$,
$(|a_1|,|a_2|,|a_3|,|b_1|,|b_2|,|b_3|)$ has to be a permutation of 
$(2,1,1,0,0,0)$ and, moreover, the terms in $F(y_i)$ 
which survives to be non-zero must have coefficients $\pm 1$.
But this is impossible, because coefficient of $y_iy_j$ is $\pm 1$ if
and only if $y_i$ and $y_j$ is connected by an edge in 
$$
\begin{array}{ccccccc}
y_1&\textrm{---}&y_2&\textrm{---}&y_3\\
&&|&&|\\
&&y_4&\textrm{---}&y_5&\textrm{---}&y_6\rlap{,}
\end{array}
$$
which does not have a cycle with three vertices. 
We therefore deduce that $|a_i|$, $|b_i|\leq 1$ when $(1/10)F\leq 1$. 

\vspace{1ex}
\statenumber\label{step-4}{\bf .}\ 
Taking some auxiliary conditions into account, we get
$$
V=\left\{
v={}^t(x_1,x_2,x_3)\,\bigg|\,
\begin{minipage}{7cm}
\setlength{\baselineskip}{.85\baselineskip}
\begin{small}
$a_i$, $b_i\in\Z$, $|a_i|$, $|b_i|\leq 1$, $b_3=0$,
$(a_1,a_2,a_3)\neq(0,0,0)$, and $F(a_i,b_j)=10$
\end{small}
\end{minipage}
\right\}.
$$
Indeed, $\tau v\in\lattice^3$ implies that $b_3$ is even (hence is
zero), and $v\not\in(\lambda/2)\Z^3$ shows
$(a_1,a_2,a_3)\neq(0,0,0)$.
Without the last condition $F(a_i,b_j)=10$, $V$ might have at most
$3^5=243$ elements. 
The rest part of the proof of Claim \ref{cla-24vectors} can be done
by searching such vectors with $F(a_i,b_j)=10$, possibly by using
computer. 

\vspace{1ex}
\statenumber\label{step-5}{\bf .}\ 
Finally, we ask for, for any $v_2\in V$, two vectors $v_1$, $v_3\in V$
such that $v_1^{\ast}Qv_2=v_2^{\ast}Qv_3=-2(\lambda+2)$ to get the matrix 
$g=(v_1,v_2,v_3)$ which is presumed to be $Q$-unitary. Due to Claim
\ref{cla-24vectors} it suffices to check the following four cases:

\vspace{1ex}
(\ref{step-5}.1)\ If $v_2={}^t(\lambda/2,1,0)$, there does not
exist such $v_1$.

\vspace{1ex}
(\ref{step-5}.2)\ If $v_2={}^t(0,0,1)$, then $v_3={}^t(0,-1,-1)$ and 
$v_1$ is either ${}^t(0,1,0)$ or ${}^t(0,-\lambda/2,-1)$; but in
this case, $(v_1,v_2,v_3)$ is not an invertible matrix.

\vspace{1ex}
(\ref{step-5}.3)\ If $v_2={}^t(0,-1,-1)$, then $v_1={}^t(0,0,1)$ and
$v_3={}^t(0,1,0)$ or ${}^t(\lambda/2,1+\lambda/2,1)$.
But, also in this case, $\det(v_1,v_2,v_3)$ is not invertible in $\Z_2$.

\vspace{1ex}
(\ref{step-5}.4)\ If $v_2={}^t(0,1,0)$, then $v_1$ is either
${}^t(1,0,0)$, ${}^t(0,-1,-1)$, or ${}^t(-\ovl{\lambda},-1,0)$, and $v_3$ is
either ${}^t(0,0,1)$ or ${}^t(-1,-1,0)$.
Among them, $v_1^{\ast}Qv_3=\lambda+2$ is satisfied only when 
$g=(v_1,v_2,v_3)=I_3$.

\vspace{1ex}
All these are checked without so much pain.
Since any member in $G_2$ is $Q$-unitary, we immediately see
that the possible $(v_1,v_2,v_3)$ are only among elements in $G_2$.
Therefore, we conclude $G_1=G_2$, as desired.
\qed

\vspace{1ex}
{\it Acknowledgments.}\ 
The first author is grateful to Professor Yves Andr\'e, to whom his
viewpoint over this work owes much.

\begin{small}

\end{small}
\end{document}